\documentclass[11pt,leqno]{article}
\usepackage{amsmath, amssymb, amsfonts, amsthm,amscd}

\swapnumbers
\theoremstyle{plain}
\newtheorem{theorem}{Theorem}[section]
\newtheorem{lemma}[theorem]{Lemma}

\newtheorem{proposition}[theorem]{Proposition}
\newtheorem{corollary}[theorem]{Corollary}

\theoremstyle{definition}
\newtheorem{definition}[theorem]{Definition}

\numberwithin{equation}{theorem}

\allowdisplaybreaks[1]
\newcommand{\leto}[1]{\stackrel{#1}{\to}}
\newcommand{\hleto}[1]{\stackrel{#1}{\hookrightarrow}}
\newcommand{\lleto}[1]{\stackrel{#1}{\leftarrow}}

\newcommand{\ds}{\displaystyle}
 
 \DeclareMathOperator{\Pic}{Pic}
\DeclareMathOperator{\Det}{Det}
 \DeclareMathOperator{\Td}{Td}
 \DeclareMathOperator{\ch}{ch}
 \newcommand{\im}{\text{Im}}
\DeclareMathOperator{\Hom}{Hom}
\begin{document}

\title{Explicit Determination of the Picard Group of Moduli Spaces of Semistable $G$-Bundles on Curves}
\author{Arzu Boysal and Shrawan Kumar}
\date{}
\maketitle


\section*{Introduction}

Let $G$ be a connected, simply-connected, simple affine algebraic group  and $\mathcal{C}_g$ be a smooth
irreducible projective curve of any genus $g\geq 1$ over $\mathbb{C}$.
Denote by $\mathfrak{M}_{\mathcal{C}_g}(G)$ the moduli space of semistable principal $G$-bundles on
$\mathcal{C}_g$.  Let $\Pic(\mathfrak{M}_{\mathcal{C}_g}(G))$ be the Picard group of
 $\mathfrak{M}_{\mathcal{C}_g}(G)$ and let $X$ be the
infinite Grassmannian
of the  affine Kac-Moody group  associated to $G$.  It is known that $\Pic(X) \simeq \mathbb{Z}$
and is generated by a homogenous line bundle $\mathfrak{L}_{\chi_0}$.  Also,
as proved by Kumar-Narasimhan  [KN], there exists a canonical
injective group homomorphism
\[
 \beta:\Pic(\mathfrak{M}_{\mathcal{C}_g}(G)) \hookrightarrow \Pic(X), \]
which takes
 $\Theta_{V}(\mathcal{C}_g,G) \mapsto\mathfrak{L}_{\chi_0}^{ m_{V}}$ for any finite dimensional representation $V$ of $G$,
where $\Theta_{V}(\mathcal{C}_g,G)$ is the theta bundle associated to the $G$-module $V$ and $m_{V}$ is its Dynkin index (cf.
Theorem 2.2). As an immediate corollary, they obtained that \[\Pic(\mathfrak{M}_{\mathcal{C}_g}(G))\simeq \mathbb{Z},\]
generalizing the corresponding result for $G=SL(n)$ proved by Drezet-Narasimhan [DN]. However, the precise image of $\beta$
was not known for non-classical $G$ excluding $G_2$. (For classical $G$ and $G_2$, see [KN], [LS], [BLS].) The main aim of
this paper is to determine the image of $\beta$ for an arbitrary $G$. It is shown that the
image of $\beta$ is generated by $\mathfrak{L}_{\chi_0}^{m_G}$, where $m_G$ is the least common multiple of the coefficients
of the coroot $\theta^\vee$  written in terms of the simple coroots, $\theta$ being the highest root of $G$ (cf. Theorem 2.4,
see also Proposition 2.3 where $m_G$ is explicitly given for each $G$).
As a consequence, we obtain that the theta bundles
$\Theta_V(\mathcal{C}_g, G)$, where $V$ runs over all the finite dimensional representations of $G$, generate $\Pic(\mathfrak{M}_{\mathcal{C}_g}(G))$ (cf. Theorem 1.3). In fact, it is shown that there is a fundamental weight $\omega_d$ such that  the theta bundle
$\Theta_{V(\omega_d)}(\mathcal{C}_g, G)$ corresponding to the irreducible highest weight $G$-module $V(\omega_d)$ with highest weight $\omega_d$ generates
$\Pic(\mathfrak{M}_{\mathcal{C}_g}(G))$ (cf. Theorem 2.4). All these fundamental weights  $\omega_d$ are explicitly determined in Proposition 2.3.

It may be mentioned that Picard group of the moduli {\it stack} of $G$-bundles is studied in [LS], [BLS], [T$_2$].

We now briefly outline the idea of the proofs.
Recall that, by a celebrated result of Narasimhan-Seshadri, the underlying real analytic space $M_g(G)$ of
$\mathfrak{M}_{\mathcal{C}_g}(G)$ admits a description as the space of representations of the fundamental group
$\pi_1(\mathcal{C}_g)$ into a fixed compact form of $G$ up to conjugation. In particular, $M_g(G)$
depends only upon $g$ and $G$ (and not on the specific choice of the projective curve $\mathcal{C}_g$).
Moreover, this description gives rise to a  standard embedding
$i_g: M_g(G) \hookrightarrow M_{g+1}(G)$.

Let $V$ be any finite dimensional representation of $G$.  We first show that the first Chern class
of the theta bundle $\Theta_V(\mathcal{C}_g, G)$ does not depend upon the choice of the smooth
projective curve $\mathcal{C}_g$, as long as $g$ is fixed (cf. Proposition~\ref{P:invcurve}).

We next show that the first Chern class of $\Theta_V(\mathcal{C}_{g+1}, G)$ restricts to the
first Chern class of $\Theta_V(\mathcal{C}_g, G)$ under the embedding $i_g$ (cf. Proposition~\ref{P:restriction}).
This result is proved by first reducing the case of general $G$ to $SL(n)$ and then reducing the
case of $SL(n)$ to $SL(2)$.  The corresponding result for $SL(2)$ is obtained by showing that the
inclusion $M_g(SL(2)) \hookrightarrow M_{g+1}(SL(2))$ induces isomorpism in cohomology
$H^2(M_{g+1}(SL(2)),\mathbb{Z}) \simeq H^2(M_g(SL(2)),\mathbb{Z})$ (cf. Proposition 1.7).
The last result  for $H^2$ with rational coefficients is fairly well known (and follows easily by observing
that the symplectic form on $M_{g+1}(G)$ restricts to the symplectic form on $M_g(G)$) but the result with
integral coefficients is more delicate and is proved in Section 4.  The proof involves the calculation
of the determinant bundle of the Poincar\'e bundle on $\mathcal{C}_g \times \mathcal{J}_{\mathcal{C}_g}$,
$\mathcal{J}_{\mathcal{C}_g}$ being the Jacobian of $\mathcal{C}_g$ which consists of the isomorphism classes
of degree $0$ line bundles on $\mathcal{C}_g$.

By  virtue of the above mentioned two propositions (Propositions~\ref{P:invcurve} and~\ref{P:restriction}),
to prove our main result determining $\Pic(\mathfrak{M}_{\mathcal{C}_g}(G))$ stated
in the first paragraph   for any $g \geq 1$, it suffices to consider the case of genus $g=1$.

In the genus $g=1$ case, $\mathfrak{M}_{\mathcal{C}_1}(G)$  admits a description as the weighted projective space $\mathbb{P}(1,a_1^\vee,a_2^\vee,\ldots , a_k^\vee),$ where $a_i^\vee$ are
the coefficients of the coroot $\theta^\vee$ written in terms of the simple coroots and $k$ is the rank of $G$ (cf. Theorems 3.1 and 3.3).  The
 ample generator of the Picard group of $\mathbb{P}(1,a_1^\vee,a_2^\vee,\ldots , a_k^\vee)$ is known to
be $\mathcal{O}_{\mathbb{P}(1,a_1^\vee,a_2^\vee,\ldots , a_k^\vee)}(m_G)$
(cf. Theorem 3.4). In
section $3$, we show that $\Theta_{V(\omega_d)}(\mathcal{C}_1,G)$ is, in fact,
$\mathcal{O}_{\mathbb{P}(1,a_1^\vee,a_2^\vee,\ldots , a_k^\vee)}(m_G)$, and hence it is the ample  generator of $\Pic(\mathfrak{M}_{\mathcal{C}_1}(G))$.
The proof makes use of the Verlinde formula determining the dimension of
the space of global
sections $H^0(\mathfrak{M}_{\mathcal{C}_g}(G), \mathfrak{L})$
(cf. Theorem 3.5).

We thank P. Belkale, L. Jeffrey and M.S. Narasimhan for some helpful correspondences/conversations. We are grateful to the referee for some very helpful comments; in particular, for pointing out Theorem 3.3, which simplified our original proof of Proposition 1.9. The second author was partially supported by the NSF.

\section{Statement of the Main Theorem and its Proof}

For a topological space $X$, $H^i(X)$ denotes the singular cohomology of $X$ with integral coefficients, unless otherwise
explicitly stated.

Let $G$ be a connected, simply-connected, simple affine algebraic group over
$\mathbb{C}$. This will be our tacit assumption on $G$ throughout the paper.
Let  $\mathcal{C}_g$ be a smooth irreducible projective curve (over $\mathbb{C}$) of genus $g$, which we
assume to be $\geq 1$.  Let $\mathfrak{M}_{\mathcal{C}_g}= \mathfrak{M}_{\mathcal{C}_g}(G)$
be the moduli space of semistable principal $G$-bundles on $\mathcal{C}_g$.

We begin by recalling the following result due to Kumar-Narasimhan [KN, Theorem 2.4]. (In loc. cit. the genus $g$ is assumed to be $\geq 2$. For the
genus $g=1$ case, the result follows from Theorems 3.1, 3.3 and 3.4.)

\begin{theorem}
With the notation as above,
\[\Pic(\mathfrak{M}_{\mathcal{C}_g}) \simeq \mathbb{Z}, \]
where
$\Pic(\mathfrak{M}_{\mathcal{C}_g})$ is the group of isomorphism classes of
algebraic line bundles on $\mathfrak{M}_{\mathcal{C}_g}$.

In particular, any nontrivial line bundle on $\mathfrak{M}_{\mathcal{C}_g}$ is ample or its inverse is ample.
\end{theorem}

\begin{definition} Let  $\mathcal{F}$ be a family of  vector bundles
on  $\mathcal{C}_g$ parametrized by a variety $X$, i.e., $\mathcal{F}$ is a  vector bundle over
 $\mathcal{C}_g \times X$. Then, the `determinant of the cohomology' gives rise to the determinant bundle $\Det(\mathcal{F})$ of the family $\mathcal{F}$, which is a line bundle over the base $X$. By definition, the fiber of   $\Det(\mathcal{F})$ over any $x\in X$ is given by the expression:
\[\Det(\mathcal{F})|_x=\wedge^{top}(H^0(\mathcal{C}_g,\mathcal{F}_x ))^* \otimes
\wedge^{top}(H^1(\mathcal{C}_g,\mathcal{F}_x)),\]
where $\mathcal{F}_x$ is the restriction of $\mathcal{F}$ to $\mathcal{C}_g\times x$  (cf., e.g.,  [L, Chap. 6, $\S$1], [KM]).

Let $\mathcal{R}(G)$ denote the set of isomorphism classes of all the finite
dimensional algebraic representations of $G$.  For any $V$ in $\mathcal{R}(G)$, we have the
$\Theta$-bundle $\Theta_V(\mathcal{C}_g) =\Theta_V(\mathcal{C}_g,G)$ on $\mathfrak{M}_{\mathcal{C}_g}$,
which is an algebraic line bundle whose fibre at any principal $G$-bundle
$E \in \mathfrak{M}_{\mathcal{C}_g}$ is given by the expression
\[ \Theta_V(\mathcal{C}_g)|_E=\wedge^{top}(H^0(\mathcal{C}_g, E_V))^* \otimes
\wedge^{top}(H^1(\mathcal{C}_g, E_V)),\] where $E_V$ is the associated vector bundle
$E\times_{G}V$ on ${\mathcal{C}_g}$. Observe that the moduli space
 $\mathfrak{M}_{\mathcal{C}_g}$ does not parametrize a universal family of $G$-bundles, however, the theta bundle  $\Theta_V(\mathcal{C}_g)$ (which is essentially the determinant bundle if there were a universal family parametrized by
 $\mathfrak{M}_{\mathcal{C}_g}$) still exists (cf. [K$_1$, $\S$3.7]).
\end{definition}

Now, we can state the main result of this paper.

\begin{theorem}\label{T:main}
\[ \Pic(\mathfrak{M}_{\mathcal{C}_g}) =<\Theta_V(\mathcal{C}_g), V \in \mathcal{R} (G)>, \]
where the notation $<\;\;>$ denotes the group generated by the elements in the bracket.
\end{theorem}

\begin{lemma}\label{L:pic=h}
\[ c:\Pic(\mathfrak{M}_{\mathcal{C}_g}) \simeq H^2(\mathfrak{M}_{\mathcal{C}_g}, \mathbb{Z}), \]
where $c$ maps any line bundle $\mathfrak{L}$ to its first Chern class $c_1(\mathfrak{L})$.

In particular,
\[H^2(\mathfrak{M}_{\mathcal{C}_g}, \mathbb{Z})\simeq  \mathbb{Z}.\]

The first Chern class of the ample generator of $\Pic(\mathfrak{M}_{\mathcal{C}_g})$ is called the
{\it positive generator} of $H^2(\mathfrak{M}_{\mathcal{C}_g}, \mathbb{Z})$.
\end{lemma}

\begin{proof}

Consider the following exact sequence of abelian groups:
\[ 0 \to \mathbb{Z} \to \mathbb{C} \leto{f} \mathbb{C^*} \to 0, \]
where $f(x)=e^{2\pi i x}$.  This gives rise to the following exact sequence of sheaves on
$\mathfrak{M}_{\mathcal{C}_g}$ endowed with the analytic topology:

\[ 0 \to \bar{\mathbb{Z}} \to \bar{\mathcal{O}}_{\mathfrak{M}_{\mathcal{C}_g}}
   \to \bar{\mathcal{O}}^*_{\mathfrak{M}_{\mathcal{C}_g}} \to 0,\]
where $\bar{\mathcal{O}}_{\mathfrak{M}_{\mathcal{C}_g}}$ is the sheaf of holomorphic
functions on $\mathfrak{M}_{\mathcal{C}_g}$, $\bar{\mathcal{O}}^*_{\mathfrak{M}_{\mathcal{C}_g}}$
is the sheaf of invertible elements of $\bar{\mathcal{O}}_{\mathfrak{M}_{\mathcal{C}_g}}$
and $\bar{\mathbb{Z}}$ is the constant  sheaf corresponding to the abelian group $\mathbb{Z}$.

The above sequence, of course,  induces the following long exact sequence in cohomology:
\[\cdots \to H^1(\mathfrak{M}_{\mathcal{C}_g}, \bar{\mathcal{O}}_{\mathfrak{M}_{\mathcal{C}_g}})
\to H^1(\mathfrak{M}_{\mathcal{C}_g}, \bar{\mathcal{O}}^*_{\mathfrak{M}_{\mathcal{C}_g}})
\leto {\bar{c}} H^2(\mathfrak{M}_{\mathcal{C}_g}, \mathbb{Z})
 \to H^2(\mathfrak{M}_{\mathcal{C}_g}, \bar{\mathcal{O}}_{\mathfrak{M}_{\mathcal{C}_g}})
\to \cdots.\]

First of all,
\[ \Pic(\mathfrak{M}_{\mathcal{C}_g})\simeq H^1(\mathfrak{M}_{\mathcal{C}_g},
\mathcal{O}^*_{\mathfrak{M}_{\mathcal{C}_g}}),\tag{1}\]
where $\mathcal{O}_{\mathfrak{M}_{\mathcal{C}_g}}$ is the sheaf of algebraic functions on
$\mathfrak{M}_{\mathcal{C}_g}$ and $\mathcal{O}^*_{\mathfrak{M}_{\mathcal{C}_g}}$ is the
subsheaf of invertible elements of $\mathcal{O}_{\mathfrak{M}_{\mathcal{C}_g}}$.

Moreover, by GAGA,   $\mathfrak{M}_{\mathcal{C}_g}$ being a projective variety,
\[ H^1(\mathfrak{M}_{\mathcal{C}_g}, \mathcal{O}^*_{\mathfrak{M}_{\mathcal{C}_g}}) \simeq
H^1(\mathfrak{M}_{\mathcal{C}_g}, \bar{\mathcal{O}}^*_{\mathfrak{M}_{\mathcal{C}_g}}),
\tag{2}\] and also, for any $p\geq 0$,
\[ H^p(\mathfrak{M}_{\mathcal{C}_g}, \mathcal{O}_{\mathfrak{M}_{\mathcal{C}_g}}) \simeq
H^p(\mathfrak{M}_{\mathcal{C}_g}, \bar{\mathcal{O}}_{\mathfrak{M}_{\mathcal{C}_g}}).\tag{3}\]

By Kumar-Narasimhan [KN, Theorem 2.8], $H^i(\mathfrak{M}_{\mathcal{C}_g},
\mathcal{O}_{\mathfrak{M}_{\mathcal{C}_g}})=0$ for $i>0$.  Hence, under the identification (1), by (2)-(3)
and the above long exact cohomology sequence,
\[\Pic (\mathfrak{M}_{\mathcal{C}_g})
\leto {{\sim}{c}} H^2(\mathfrak{M}_{\mathcal{C}_g}, \mathbb{Z}),\]
where $c$ is the map $\bar{c}$ under the above identifications. Moreover,  as is well known, $c$
is the first Chern class map.
\end{proof}

\vspace{5pt}

Let us fix a  maximal compact subgroup $K$ of $G$.  Denote the Riemann surface with $g$ handles,
considered only as  a topological manifold,   by $C_g$.  Thus, the underlying topological manifold
of $\mathcal{C}_g$ is $C_g$.
Define $M_g(G):=\varphi^{-1}(1)/\mbox{Ad}K$, where $\varphi :K^{2g}\rightarrow K$
is the commutator map $\varphi (k_1,k_2,\ldots, k_{2g}) = [k_1,k_2][k_3,k_4]\cdots[k_{2g-1},k_{2g}]$
and $\varphi^{-1}(1)/\mbox{Ad}K$ refers to the quotient of $\varphi^{-1}(1)$ by $K$ under the diagonal
adjoint action of $K$ on $K^{2g}$.

Now, we recall the following fundamental result due to Narasimhan-Seshadri [NS] for vector bundles and
extended for an arbitrary $G$ by Ramanathan [R$_1$, R$_2$].

Consider the standard generators $a_1,b_1,a_2,b_2, \ldots, a_g,b_g$ of
$\pi_1(\mathcal{C}_g)$  (cf. [N, $\S14$]). Then, we have the presentation:
$$\pi_1(\mathcal{C}_g)=F[a_1,\ldots, a_g, b_1,\ldots, b_g]/<[a_1,b_1]\cdots [a_g,b_g]>,$$ where $
F[a_1,\ldots, a_g, b_1,\ldots, b_g]$ denotes the free group generated by
$a_1,\ldots, a_g$, $
b_1,\ldots, b_g$ and $<[a_1,b_1]\cdots [a_g,b_g]>$ denotes the normal subgroup generated by the single element $[a_1,b_1]\cdots [a_g,b_g]$.

\begin{theorem}
Having chosen the standard generators $a_1,b_1,a_2,b_2, \ldots, a_g,b_g$ of
$\pi_1(\mathcal{C}_g)$, there exists a canonical  isomorphism of real analytic spaces:
\[ \theta_{\mathcal{C}_g} (G): M_g(G) \simeq \mathfrak{M}_{\mathcal{C}_g} (G). \]
\end{theorem}

In the sequel, we will often make this identification.

\begin{proposition}\label{P:invcurve}
For any $V \in \mathcal{R}(G)$, $c(\Theta_V(\mathcal{C}_g,G))$, under the above identification
$\theta_{\mathcal{C}_g} (G)$,  does not depend on the choice of the projective variety structure
$\mathcal{C}_g$ on the Riemann surface $C_g$ for any fixed $g$.
\end{proposition}

\begin{proof}
Let $\rho:G \to SL(V)$ be the given representation.  By taking a $K$-invariant Hermitian form on $V$ we
get $\rho(K) \subset SU(n)$, where $n=\mbox{dim}\,V$. For any principal $G$-bundle $E$ on $\mathcal{C}_g$,
let $E_{SL(V)}$  be the  principal $SL(V)$-bundle over $\mathcal{C}_g$ obtained by the extension of the
structure group via $\rho$. Then, if $E$ is semistable, so is $E_{SL(V)}$,  giving rise to a variety
morphism $ \hat{\rho}: \mathfrak{M}_{\mathcal{C}_g}(G) \to \mathfrak{M}_{\mathcal{C}_g}(SL(V))$ (cf. [RR, Theorem 3.18]).
Hence, we get the commutative diagram:
 \[ \begin{aligned}    &\mathfrak{M}_{\mathcal{C}_g}(G) \leto{\hat{\rho}}
\mathfrak{M}_{\mathcal{C}_g}(SL(V))\\
& \hspace{5pt} \uparrow \hspace{70pt}
        \uparrow  \\
     & M_g(G)  \leto{\bar{\rho}} M_g(SL(V)),
\end{aligned} \tag{D$_1$}\]
where $\bar{\rho}$ is induced from the commutative diagram:
 \[  \begin{CD}
       K^{2g} @>{\varphi} >> K \\
       @VV{\rho^{\times 2g}}V    @VV{\rho}V \\
       SU(n)^{2g} @>>\varphi>  SU(n).
    \end{CD} \]
The diagram (D$_1$)  induces the following commutative diagram in cohomology:
 \[  \begin{CD}
      H^2(\mathfrak{M}_{\mathcal{C}_g}(SL(V)), \mathbb{Z}) @> {\hat{\rho}^*}
       >>  H^2(\mathfrak{M}_{\mathcal{C}_g}(G), \mathbb{Z}) \\
      @VV\|V    @VV\|V \\
      H^2(M_g(SL(V)), \mathbb{Z}) @>>{\bar{\rho}^*}>   H^2(M_g(G), \mathbb{Z}).
     \end{CD}   \tag{D$_2$}\]

By the construction of the $\Theta$-bundle, $\hat{\rho}^*(\Theta_V(\mathcal{C}_g,SL(V)))
=\Theta_V(\mathcal{C}_g, G),$ where $\hat{\rho}^*$ also denotes the pullback of  line bundles
and $V$ is thought of as the standard representation of $ SL(V)$.

Thus, using the functoriality of the Chern class, we get
\[ \hat{\rho}^*(c(\Theta_V(\mathcal{C}_g, SL(V))))
  =c(\Theta_V(\mathcal{C}_g, G)). \tag{1}\]
By Drezet-Narasimhan [DN], $c(\Theta_V(\mathcal{C}_g, SL(V)))$ is the unique positive generator of
$H^2(\mathfrak{M}_{\mathcal{C}_g}(SL(V)), \mathbb{Z})$ and thus is independent of the choice of $\mathcal{C}_g$
under the identification $\theta_{\mathcal{C}_g} (SL(V))^*$.  Consequently, by (1) and the above
commutative diagram (D$_2$), $c(\Theta_V(\mathcal{C}_g, G))$ is independent of the choice of $\mathcal{C}_g$.
\end{proof}

From now on we will denote the cohomology class $c(\Theta_V(\mathcal{C}_g,G))$ in $H^2(M_g(G), \mathbb{Z})$,
under the identification $\theta_{\mathcal{C}_g} (G)^*$, by $c(\Theta_V(g,G))$.

Consider the embedding \[ i_g=i_g(G): M_g(G) \hookrightarrow M_{g+1}(G) \]
induced by the inclusion of $K^{2g} \to K^{2g+2}$ via $(k_1,\ldots, k_{2g}) \mapsto
(k_1,\ldots,k_{2g},1,1)$.

By virtue of the map $i_g$, we will identify $M_g(G)$ as a subspace of $M_{g+1}(G)$.
In particular, we get the following induced sequence of maps in the second cohomology.
\[    H^2(M_1(G), \mathbb{Z}) \lleto{i_1^*}  H^2(M_2(G), \mathbb{Z})
      \lleto{i_2^*}  H^2(M_3(G), \mathbb{Z}) \lleto{i_3^*} \cdots \,\cdot \]

\begin{proposition}\label{P:sl2case}
For $G=SL(2)$, the maps $i_g^*: H^2( M_{g+1}(G), \mathbb{Z})
\to H^2( M_g (G), \mathbb{Z})$ are isomorphisms for any $g\geq 1$.

In particular,  $i_g^*$ takes the positive generator of $ H^2( M_{g+1}(SL(2)), \mathbb{Z})$ to  the positive
generator of $ H^2( M_{g}(SL(2)), \mathbb{Z}).$
\end{proposition}

We shall prove this proposition in Section 4.

\begin{proposition}\label{P:restriction} For any $V\in \mathcal{R}(G)$  and any $g\geq 1$,
$i_g^*(c(\Theta_V(g+1,G)))=c(\Theta_V(g,G))$.
\end{proposition}

\begin{proof}
We first claim that it suffices to prove the above proposition for $G=SL(n)$ and the
standard $n$-dimensional representation $V$ of $SL(n)$.

Let $\rho: G \to SL(V)$ be the given representation.
Consider the following commutative diagram:
\[ \begin{aligned}
   &M_g(G) & \hleto{i_g} &M_{g+1}(G)\\
   &\bar{\rho} \downarrow & & \hspace{5pt} \downarrow \bar{\rho} \\
   &M_g(SL(V)) &\hleto{i_g} &M_{g+1}(SL(V)),
\end{aligned} \]
where $ \bar{\rho}$ is the map defined in the proof of Proposition~\ref{P:invcurve}.
It induces the commutative diagram:
\[ H^2(M_g(G),\mathbb{Z}) \lleto{i_g^*} H^2(M_{g+1}(G),\mathbb{Z}) \]
\[ {\bar{\rho}}^* \uparrow \hspace{60pt} {\bar{\rho}}^* \uparrow \]
\[ H^2(M_g(SL(V)),\mathbb{Z}) \lleto{i_g^*} H^2(M_{g+1}(SL(V)),\mathbb{Z}).\]
Therefore, using the commutativity of the above diagram and equation (1) of
Proposition~\ref{P:invcurve},
supposing that $i_g^*(c(\Theta_V(g+1,SL(V))))=c(\Theta_V(g,SL(V)))$, we get
$i_g^*(c(\Theta_V(g+1,G)))=c(\Theta_V(g,G))$.  Hence, Proposition~\ref{P:restriction}
is established for any $G$ provided we assume its validity for
$G=SL(V)$ and its standard representation in $V$.

We further reduce the proposition from $SL(n)$ to  $SL(2)$. As in the proof of Proposition~\ref{P:invcurve},
consider the mappings
\[ \bar{\rho}: M_g(SL(2)) \rightarrow M_{g}(SL(n)), \,\,\text{and}\]
\[\hat{\rho}: \mathfrak{M}_{\mathcal{C}_g}(SL(2)) \rightarrow
\mathfrak{M}_{\mathcal{C}_g}(SL(n))\]
induced by the inclusions
\[
 SU(2) \to SU(n) \,\,\,\text{and}\, SL(2) \to SL(n),
   \]
given by $m  \mapsto diag(m, 1, \ldots, 1).$

The maps $\bar{\rho}$ and $\hat{\rho}$ induce the commutative diagram:
 \[  \begin{aligned}    &
 H^2(M_g(SL(n)),\mathbb{Z})  \leto{{\bar{\rho}}^*}
  H^2(M_{g}(SL(2)),\mathbb{Z})  \\
& \hspace{35pt} \uparrow \hspace{115pt}
        \uparrow  \\
 & H^2(\mathfrak{M}_{\mathcal{C}_g}(SL(n)),\mathbb{Z}) \leto{{\hat{\rho}}^*}
 H^2(\mathfrak{M}_{\mathcal{C}_g}(SL(2)),\mathbb{Z}).
\end{aligned}\]
By the construction of the $\Theta$-bundle, $\hat{\rho}^*\bigl(\Theta_V(\mathcal{C}_g,
SL(n))\bigr) =\Theta_{V_2}(\mathcal{C}_g, SL(2))$, where $V_2$  is the standard $2$-dimensional representation of $SL(2)$.

Thus, using the functoriality of the Chern class, we get
\[ {\hat{\rho}}^*(c(\Theta_V(\mathcal{C}_g, SL(n))))=c(\Theta_{V_2}(\mathcal{C}_g, SL(2))).\tag{1}\]

Using one more time the result of Drezet-Narasimhan that $c(\Theta_V(\mathcal{C}_g, SL(n)))$ is the unique
positive generator of $H^2(\mathfrak{M}_{\mathcal{C}_g}(SL(n)))$ for any $n$ (cf. Proof of
Proposition~\ref{P:invcurve}), we see that $ {\hat{\rho}}^*$ is surjective and hence an isomorphism by
Lemma~\ref{L:pic=h}.

Consider the following commutative diagram:
\[ H^2(M_g(SL(n)),\mathbb{Z}) \lleto{{i_g^*}} H^2(M_{g+1}(SL(n)),\mathbb{Z})\]
\[ {\bar{\rho}}^* \downarrow \hspace{60pt} {\bar{\rho}}^* \downarrow \]
\[ H^2(M_g(SL(2)),\mathbb{Z}) \lleto{i_g^*} H^2(M_{g+1}(SL(2)),\mathbb{Z}).\]
Suppose that the proposition is true for $G=SL(2)$ and the standard representation
$V_2$, i.e.,
\[i_g^*(c(\Theta_{V_2}(g+1,SL(2))))=c(\Theta_{V_2}(g,SL(2))). \tag{2}\]
  Then, using the
commutativity of the above diagram and (1) together with  the fact
that ${\bar{\rho}}^*$ is an isomorphism, we get that
$i_g^*(c(\Theta_V(g+1,SL(n))))=c(\Theta_V(g,SL(n)))$. Finally, (2)
follows from the result of Drezet-Narasimhan cited above and
Proposition~\ref{P:sl2case}. Hence the proposition is established
for any $G$ (once we prove Proposition 1.7).
\end{proof}

\begin{proposition}\label{P:ellcase}
For $g=1$, Theorem~\ref{T:main} is true.
\end{proposition}

The proof of this proposition will be given in Section 3.

\begin{proof}[Proof of Theorem~\ref{T:main}]
Denote the subgroup $<\Theta_V(\mathcal{C}_g,G), V \in \mathcal{R}(G)>$ of
$\Pic(\mathfrak{M}_{\mathcal{C}_g}(G))$
by
$\Pic^{\Theta}(\mathfrak{M}_{\mathcal{C}_g}(G))$.

Set $H^2_{\Theta}(M_g(G)) :=
c(\Pic^{\Theta}(\mathfrak{M}_{\mathcal{C}_g}(G)))$.  By virtue of
Proposition~\ref{P:invcurve}, this is well defined, i.e.,
$H^2_{\Theta}(M_g(G))$ does not depend upon the choice of the
projective variety structure $\mathcal{C}_g$ on $C_g$.  Moreover,
by Proposition~\ref{P:restriction},
$i_g^*(H^2_{\Theta}(M_{g+1}(G)))=H^2_{\Theta}(M_{g}(G))$.

Thus, we get the following commutative diagram, where the upward arrows are inclusions and the maps in the bottom horizontal sequence are induced from the maps $i_g^*$.

\[ \begin{aligned}
      & H^2(M_1(G)) \lleto{i_1^*} H^2(M_2(G))
                    \lleto{i_2^*} H^2(M_3(G)) \lleto{i_3^*} \cdots \\
      & \hspace{20pt} \uparrow  \hspace{60pt} \uparrow \hspace{60pt}
        \uparrow  \\
      & H^2_{\Theta}(M_1(G)) \twoheadleftarrow H^2_{\Theta}(M_2(G)) \twoheadleftarrow
        H^2_{\Theta}(M_3(G)) \twoheadleftarrow \cdots .
 \end{aligned} \]

By Proposition~\ref{P:ellcase} and Lemma~\ref{L:pic=h}, $H^2(M_1(G))=H^2_{\Theta}(M_1(G))$.
Then, $i_1^*$ is surjective and hence an isomorphism (by using Lemma~\ref{L:pic=h} again).
Thus, by the commutativity of the above diagram, the inclusion  $H^2_{\Theta}(M_2(G))
\hookrightarrow H^2(M_2(G))$ is an isomorphism. Arguing the same way, we get that  $H^2(M_g(G))$
$=
H^2_{\Theta}(M_g(G))$ for all $g$. This completes the proof of the theorem by virtue of the isomorphism
$c$ of Lemma~\ref{L:pic=h}.
\end{proof}


\section{Comparison of the Picard Groups of $\mathfrak{M}_{\mathcal{C}_g}$ and the Infinite Grassmannian}

As earlier, let $G$ be a connected, simply-connected, simple affine algebraic group over $\mathbb{C}$. We fix a
Borel subgroup $B$ of $G$ and a maximal torus $T\subset B$. Let $\mathfrak h$ (resp. $\mathfrak b$) be the Lie
algebra of $T$ (resp. $B$). Let $\Delta^+ \subset \mathfrak h^*$ be the set of positive roots (i.e., the roots
of $\mathfrak b$  with respect to $\mathfrak h$) and let $\{\omega_i\}_{1\leq i\leq k} \subset \mathfrak h^*$ be
the set of fundamental weights, where $k$ is the rank of $G$. As earlier, $\mathcal{R}(G)$ denotes the set of
isomorphism classes of all the finite dimensional algebraic representations of $G$. This is a semigroup under
the direct sum of two representations. Let $R(G)$ denote the associated Grothendieck group. Then, $R(G)$ is a
ring, where the product is induced from the tensor product of two representations. Then, the  fundamental
representations $\{ V(\omega_i)\}_{1 \leq i \leq k}$ generate the representation ring $R(G)$ as a ring [A].

Let $X$ be the infinite Grassmannian  associated to the affine Kac-Moody group
$\mathcal G$ corresponding to $G$, i.e., $X:= \mathcal G/\mathcal P$, where
$\mathcal P$ is the standard maximal parabolic subgroup of $\mathcal G$
 (cf. [K$_2$, $\S$13.2.12]; in loc. cit., $X$ is denoted by $\mathcal Y=\mathcal X^Y$).  It is known that
 $\Pic(X)$ is isomorphic to $\mathbb{Z}$ and is generated by the homogenous line bundle
 $\mathfrak{L}_{\chi_0}$ (cf. [K$_2$, Proposition 13.2.19]).

We recall the following definition from [D,$\S$2].

\begin{definition}    Let $\mathfrak{g}_1$ and
 $\mathfrak{g}_2$ be two (finite dimensional) complex simple
 Lie algebras and $\varphi : \mathfrak{g}_1 \to \mathfrak{g}_2$
 be a Lie algebra  homomorphism.  There exists a unique
 number $m_\varphi \in \mathbb C$, called the {\it Dynkin index}
 of the homomorphism $\varphi$, satisfying
 $$
 \langle \varphi (x),\varphi (y) \rangle
 = m_\varphi \langle x,y \rangle, \text{ for all }x,y \in
 \mathfrak{g}_1,
 $$
 where $\langle , \rangle$ is the Killing form on
 $\mathfrak{g}_1$ (and $\mathfrak{g}_2$) normalized so that
 $\langle \theta,\theta \rangle=2$ for the highest root
 $\theta$.

For a Lie algebra $\mathfrak g_1$ as above and a
finite dimensional representation $V$ of $\mathfrak g_1$, by the {\it Dynkin
index} $m_V$ of $V$, we mean the   Dynkin
index of the
Lie algebra homomorphism $\rho : \mathfrak{g}_1 \to sl(V)$, where  $sl(V)$
is the Lie algebra of trace $0$ endomorphisms of $V$.

Then, for any two  finite
 dimensional representations  $V$ and $W$ of $\mathfrak{g}_1$, we have,
by [D, Chap. 1, $\S2$] or [KN, Lemma 4.5],
 \[
 m_{V\otimes W} = m_V\dim W + m_W \dim V.
\tag{1} \]
\end{definition}

We recall the following main result of Kumar-Narasimhan [KN, Theorem 2.4].

\begin{theorem} \label{T:affineX} There exists a `natural'  injective group homomorphism
\[
 \beta: \Pic(\mathfrak{M}_{\mathcal{C}_g}(G)) \hookrightarrow \Pic(X). \]
Moreover, by [KNR, Theorem 5.4] (see also [Fa]), for any $V\in \mathcal R(G)$,
 \[\beta(\Theta_{V}(\mathcal{C}_g,G)) =\mathfrak{L}_{\chi_0}^{\otimes m_{V}},
\tag{1}  \]
where $V$ is thought of as a module for $\mathfrak g$ under differentiation and $m_{V}$ is its Dynkin index.
\end{theorem}

We also recall the following result from [D, Table 5],  [KN, Proposition 4.7], or [LS, $\S2$].

\begin{proposition}\label{P:dindex}  For any simple Lie algebra $\mathfrak{g}$, there exists a
(not unique in general) fundamental weight  $\omega_d$ such that
$m_{V(\omega_d)}$ divides each of $\{m_{V(\omega_i)}\}_{1\leq i \leq k}$. Thus, by (1) of Definition 2.1, $
m_{V(\omega_d)}$ divides $m_V$ for any $V\in \mathcal{R}(G)$.

The following table  gives the list of all such $\omega_d$'s and the corresponding Dynkin index
$m_{V(\omega_d)}$.

\[ \begin{array}{ccc}
\mbox{Type of G}                     & \omega_d              & m_{V(\omega_d)}\\
A_k \;(k\geq 1)                      & \omega_1, \omega_k              & 1 \\
C_k \;(k\geq 2)                      & \omega_1              & 1 \\

B_k \;(k\geq 3)                      & \omega_1              & 2 \\
D_k \;(k\geq 4)                      & \omega_1              & 2 \\
G_2                                  & \omega_1              & 2 \\
F_4                                  & \omega_4              & 6 \\
E_6                                  & \omega_1, \omega_6 & 6 \\
E_7                                  & \omega_7              & 12 \\
E_8                                  & \omega_8              & 60.
\end{array} \]
For $B_3$, $\omega_3$ also satisfies $m_{V(\omega_3)}=2$; for $D_4$,
$\omega_3$ and $\omega_4$ both have  $m_{V(\omega_3)}=m_{V(\omega_4)}=2$.
\end{proposition}

Let $\theta$ be the highest root of $G$. Observe that, for any $G$, $m_{V(\omega_d)}$ is the least common
multiple of the coefficients of the  coroot $\theta^\vee$  written in terms of the simple coroots. We shall
denote $m_{V(\omega_d)}$ by $m_G$.

Combining the above result with Theorem~\ref{T:main}, we get the following.

\begin{theorem} For any $\mathcal{C}_g$ with $g\geq 1$ and $G$ as in Section 1, the Picard group
$ \Pic(\mathfrak{M}_{\mathcal{C}_g}(G))$ is freely generated by the $\Theta$-bundle
$\Theta_{V(\omega_d)}(\mathcal{C}_g,G)$, where $\omega_d$ is any fundamental weight as in the above
proposition.

In particular,
\[\im (\beta) \,\text{is freely generated by }\, \mathfrak{L}_{\chi_0}^{\otimes m_G}. \tag{1}\]
\end{theorem}

\begin{proof} By Theorem~\ref{T:main},
\[\Pic(\mathfrak{M}_{\mathcal{C}_g}(G))=<\Theta_V(\mathcal{C}_g, G), V\in
\mathcal{R}(G)>.\]
Thus, by Theorem~\ref{T:affineX} and Proposition 2.3,
\[\im (\beta)=\,<\mathfrak{L}_{\chi_0}^{\otimes m_{V}}, V\in \mathcal{R}(G)>
\,= \,<\mathfrak{L}_{\chi_0}^{\otimes m_G}>
.\]
This proves (1).

Since $\beta$ is injective, by the above description of $\im (\beta)$,
$\Theta_{V(\omega_d)}(\mathcal{C}_g, G)$ freely generates
$\Pic(\mathfrak{M}_{\mathcal{C}_g}(G))$, proving the theorem.
\end{proof}

Following the same argument as in [So, $\S$4], using Theorem 2.4 and Proposition 2.3, we get the following corollary for genus $g\geq 2$. For genus $g=1$,
use Theorems 3.1 and 3.3 together with [BR, Theorem 7.1.d]. This corollary is due to [BLS], [So].

\begin{corollary} Let $G$ be any group  and $\mathcal{C}_g$ be any curve as in Section 1. Then, the moduli space $\mathfrak{M}_{\mathcal{C}_g}(G)$ is locally factorial if and only if $G$ is of type $A_k$ ($k\geq 1$) or  $C_k$ ($k\geq 2$).
\end{corollary}

\section{Proof of Proposition~\ref{P:ellcase}}

Let $G$ be as in the beginning of Section 1. In this section, we  identify
 $\mathfrak{M}_{\mathcal{C}_1}(G)$ with a weighted
projective space and show that the generator of $\Pic(\mathfrak{M}_{\mathcal{C}_1}(G))$ is
$\Theta_{V(\omega_d)}(\mathcal{C}_1, G)$ as claimed.

  We recall the  following theorem  due independently to Laszlo [La, Theorem 4.16] and
Friedman-Morgan-Witten [FMW, $\S$2].

\begin{theorem}
Let $\mathcal{C}_1$
be a smooth, irreducible projective curve of genus $1$. Then, there is a
natural variety isomorphim   between the moduli space
$\mathfrak{M}_{\mathcal{C}_1}(G)$ and $(\mathcal{C}_1
\otimes_{\mathbb{Z}}Q^\vee)/W$, where $Q^\vee$ is the coroot lattice
of $G$ and $W$ is its Weyl group acting canonically on $Q^\vee$ (and acting trivially on $\mathcal{C}_1$).
\end{theorem}

\begin{definition}
Let $N=(n_0,\ldots,n_k)$ be a $k+1$-tuple of positive integers.  Consider the polynomial ring $\Bbb C[z_0,\ldots,z_k]$
graded by $\mbox{deg}\,z_i=n_i$.  The scheme $\mbox{Proj}(\Bbb C[z_0,\ldots,z_k])$ is said to be the weighted
projective space of type $N$ and we  denote it by $\mathbb{P}(N)$.

Consider the standard (nonweighted) projective space $\mathbb{P}^k:=$
$\mbox{Proj}(\Bbb C[w_0,\ldots,w_k])$, where each $\mbox{deg}\,w_i=1$. Then, the graded algebra homomorphism $\Bbb C[z_0,\ldots,z_k] \to \Bbb C[w_0,\ldots,w_k],
z_i\mapsto w_i^{n_i}$, induces a morphism $\delta: \mathbb{P}^k \to \mathbb{P}(N)$.
\end{definition}

The following theorem is due to Looijenga [Lo]. His proof had a gap; a complete proof of a more general result is outlined by Bernshtein-Shvartsman [BSh].

\begin{theorem}\label{T:wps}
Let $\mathcal{C}_1$ be an elliptic curve.  Then, the variety
$(\mathcal{C}_1\otimes_{ \mathbb{Z}}Q^\vee)/W$ is the weighted projective space of type
$(1,a_1^\vee,a_2^\vee,\ldots , a_k^\vee)$, where $a_i^\vee$ are the coefficients of the coroot $\theta^\vee$ written in
terms of the simple coroots $\{\alpha_i^\vee\}$ (and, as earlier,  $k$ is the rank of $G$).
\end{theorem}

The following table lists the weighted projective space isomorphic to
$\mathfrak{M}_{\mathcal{C}_1}(G)$ corresponding to any $G$.
In this table the entries beyond $1$ are precisely the numbers $(a_1^\vee,a_2^\vee,\ldots , a_k^\vee)$ following the convention as in Bourbaki [B, Planche I-IX].
\[ \begin{array}{cc}
\mbox{Type of G}                              & \mbox{Type of the weighted projective space}\\
A_k \;(k\geq 1), \;\;C_k \;(k\geq 2)          & (1,1,1,\ldots,1) \\
B_k  \;(k\geq 3)                    & (1,1,2,\ldots,2,1) \\
D_k  \;(k\geq 4)                   & (1,1,2,\ldots,2,1,1)\\
G_2                      & (1,1,2) \\
F_4                      & (1,2,3,2,1) \\
E_6                      & (1,1,2,2,3,2,1)\\
E_7                      & (1,2,2,3,4,3,2,1) \\
E_8                      & (1,2,3,4,6,5,4,3,2) .
\end{array} \]

We recall the following result from the theory of weighted projective spaces (see, e.g.,  Beltrametti-Robbiano [BR, Lemma 3B.2.c and Theorem 7.1.c]).

\begin{theorem}\label{T:picwps}
Let $N=(n_0,\ldots, n_k)$ and assume $\mbox{gcd}\{ n_0,\ldots,n_k \}=1$. Then, we have the following.

(a)  $\Pic(\mathbb{P}(N))\simeq\mathbb{Z}$. In fact,
the morphism $\delta$ of Definition 3.2 induces an injective map
$\delta^*: \Pic(\mathbb{P}(N)) \to \Pic(\mathbb{P}^k)$.

Moreover, the ample generator of $\Pic(\mathbb{P}(N))$ maps to
$\mathcal{O}_{\mathbb{P}^k}(s)$ under $\delta^*$, where $s$ is the least common multiple of $\{ n_0,\ldots,n_k \}$. We denote this ample generator by
$\mathcal{O}_{\mathbb{P}(N)}(s)$.

(b) For any $d\geq 0$,
$$H^0(\mathbb{P}(N),
\mathcal{O}_{\mathbb{P}(N)}(s)^{\otimes d})=
\Bbb C[z_0,\ldots,z_k]_{ds},$$
where $\Bbb C[z_0,\ldots,z_k]_{ds}$ denotes the subspace of $\Bbb C[z_0,\ldots,z_k]$ consisting of homogeneous elements of degree $ds$.
\end{theorem}

Using Theorems 3.1, \ref{T:wps} and \ref{T:picwps} and the fact that the least common multiple of the numbers $\{1,a_1^\vee,a_2^\vee,\ldots ,
a_k^\vee \}$ for each $G$ is the Dynkin index $m_G=m_{V(\omega_d)}$, we have
\[\Theta_{V(\omega_d)}(\mathcal{C}_1, G)
=\mathcal{O}_{\mathbb{P}
(1,a_1^\vee,a_2^\vee,\ldots , a_k^\vee)}(m_G)^{\otimes p}\tag{$*$} \]
 for some positive integer $p$.
The value of $m_G$ is given in Proposition~\ref{P:dindex} for any $G$.

We recall the following basic result, the first part of which is due independently to Beauville-Laszlo [BL],
Faltings [Fa] and Kumar-Narasimhan-Ramanathan [KNR]. The second part of the theorem as in (1) is the celebrated
Verlinde formula for the dimension of the space of conformal blocks essentially due to Tsuchiya-Ueno-Yamada
[TUY] (together with works [Fa, Appendix] and [T$_1$]).

\begin{theorem}\label{T:cbv}
For any ample line bundle $\mathfrak{L}\in
\Pic(\mathfrak{M}_{\mathcal{C}_g}(G))$  and $\ell\geq 0$, there is an
isomorphism (canonical up to scalar multiples):
\[ H^0(\mathfrak{M}_{\mathcal{C}_g}(G), \mathfrak L^{\otimes \ell})
\simeq L(\mathcal{C}_g,\ell m_{\mathfrak L}), \] where $L(\mathcal{C}_g,\ell)$ is the space of conformal blocks
corresponding to the one marked point on $\mathcal{C}_g$ and trivial representation attached to it with central
charge $\ell$ (cf., e.g., [TUY] for the definition of conformal blocks) and $m_{\mathfrak L}$ is the positive
integer such that $\beta (\mathfrak L)= \mathfrak L_{\chi_0}^{\otimes m_{\mathfrak L}}$, $\beta$  being the map
as in Theorem~\ref{T:affineX}.

Moreover, the dimension $ F_g(\ell)$ of the space $L(\mathcal{C}_g,\ell)$ is given by the following Verlinde
formula:

\[ F_g(\ell)=t_\ell^{g-1} \ds \sum_{\mu \in P_{\ell}} \ds \prod_{\alpha \in
\Delta_{+}} \vert 2 \sin(\frac{\pi}{\ell +h}<\alpha,\mu+\rho>)\vert^{2-2g}, \tag{1}\] where
\[ \begin{array}{l}
 <\;,\;>:= \mbox{Killing form on $\mathfrak h^*$ normalized so that $<\theta,\theta>=2$ for the highest root $\theta$} \\
 \Delta_+ := \mbox{the set of positive roots}, \\
 P_{\ell}:= \{ \mbox{dominant integral weights} \;\;\mu | <\mu, {\theta}>\leq \ell \}, \\
 \rho := \mbox{ half sum of positive roots}, \\
 h:= <\rho, {\theta}>+1, \;  \mbox{the dual Coxeter number}, \\
 t_\ell:= (\ell +h)^{\mbox{rank}\, G}(\#P/Q_{lg}),
\end{array} \]
and $P$ is the weight lattice and $Q_{lg}$ is the sublatttice of the root lattice $Q$ generated by the long
roots.
\end{theorem}

In fact, we only need to use the above theorem  for the case of genus $g=1$. For $g=1$, the Verlinde formula
(1) clearly reduces to the identity:
\[F_1(\ell)= \# P_{\ell}.\]
Of course,
$$P_\ell=\{(n_1, \ldots, n_k)\in (\Bbb Z_+)^k: \sum_{i=1}^k\,n_ia_i^\vee\leq \ell\}.$$

\begin{proof}[Proof of Proposition 1.9] Using the specialization of Theorem~\ref{T:cbv}
to $g=1$, we see that $$ \mbox{dim}\;
H^0(\mathfrak{M}_{\mathcal{C}_1}(G),\Theta_{V(\omega_d)}(\mathcal{C}_1, G))=\# P_{m_G}.$$
On the other hand, by Theorems~ 3.1, \ref{T:wps} and ~\ref{T:picwps}(b), $$ \mbox{dim}\;
H^0(\mathfrak{M}_{\mathcal{C}_1}(G), \mathcal{O}_{\mathbb{P}(1,a_1^\vee,a_2^\vee,\ldots ,
a_k^\vee)}(m_G)^{\otimes p}) = \mbox{dim}({\Bbb C}[z_0,\ldots,z_k]_{pm_G})= \# P_{pm_G}.
$$
Hence, in the equation ($*$) following Theorem~\ref{T:picwps}, $p=1$ and
$\Theta_{V(\omega_d)}(\mathcal{C}_1,G)$ is the (ample)  generator of
$\Pic(\mathfrak{M}_{\mathcal{C}_1}(G))$. This proves Proposition 1.9.  \end{proof}


\section{Proof of Proposition 1.7}

In this section, we take $G=SL(2)$ and abbreviate $\mathfrak{M}_{\mathcal{C}_g}(SL(2))$ by $\mathfrak{M}_{\mathcal{C}_g}$ etc.
Let $\mathfrak{M}^{\mbox{\tiny red}}_{\mathcal{C}_g}$ be the closed subvariety of the moduli space $\mathfrak{M}_{\mathcal{C}_g}$
consisting of decomposable bundles on $\mathcal{C}_g$ (which are semistable of rank-$2$ with trivial determinant).  Let
$\mathfrak{J}_{\mathcal{C}_g}$ be the Jacobian of $\mathcal{C}_g$. Recall that the underlying set of the variety
$\mathfrak{J}_{\mathcal{C}_g}$ consists of all the isomorphism classes of line bundles on $\mathcal{C}_g$ with trivial
first Chern class. Then, there is a surjective morphism
$\xi=\xi_{\mathcal{C}_g}: \mathfrak{J}_{\mathcal{C}_g} \to \mathfrak{M}^{\mbox{\tiny red}}_{\mathcal{C}_g}\subset \mathfrak{M}_{\mathcal{C}_g}$,  taking $\mathfrak L \mapsto
\mathfrak L \oplus \mathfrak L^{-1}$. Moreover, $\xi^{-1}(\xi(\mathfrak L))=
\{ \mathfrak L, \mathfrak L^{-1}\}$. The Jacobian
$\mathfrak{J}_{\mathcal{C}_g}$ admits the involution $\tau$ taking
$\mathfrak L \mapsto \mathfrak L^{-1}$.

Let $T$ be a maximal torus of the maximal compact subgroup $SU(2)$ of $SL(2)$, which we take to be the diagonal subgroup of  $SU(2)$. Similar to the identification $\theta_{\mathcal{C}_g}$ as in Theorem
1.5, setting $J_g:=T^{2g}$, there is an isomorphism of real analytic spaces
$\bar{\theta}_{\mathcal{C}_g}: J_g \to \mathfrak{J}_{\mathcal{C}_g}$ making the following diagram commutative:
\[ J_g  \leto{\bar{\theta}_{\mathcal{C}_g}} \mathfrak{J}_{\mathcal{C}_g} \]
\[f_g \downarrow  \hspace{30pt} \downarrow \xi_{\mathcal{C}_g} \tag{E}\]
\[ M_g\leto{{\theta}_{\mathcal{C}_g}} \mathfrak{M}_{\mathcal{C}_g},\]
where $f_g: J_g \to M_g$ is induced from the standard inclusion
$T^{2g} \subset SU(2)^{2g}.$ We will explicitly describe the isomorphism
$\bar{\theta}_{\mathcal{C}_g}$ in the proof of the following lemma.

Recall the definition of the map $i_g:M_g \to M_{g+1}$ from Section 1 and let
$r_g: J_g \to J_{g+1}$ be the map $(t_1, \ldots, t_{2g})\mapsto
(t_1, \ldots, t_{2g},1,1).$ Then, we have the following commutative diagram:
\[ \begin{aligned}    & J_g  \leto{f_g}\,\,\,\,\,\,\, M_{g} \\
 r_g &\downarrow  \hspace{40pt} \downarrow i_g \\
&J_{g+1} \leto{f_{g+1}} M_{g+1}.
\end{aligned} \tag{F}\]
Let $x_{g+1}$ denote the positive generator of $H^2(M_{g+1}, \mathbb{Z})$.  Then, by Lemma 1.4 and Theorem 1.5,
\[i_g^*(x_{g+1})=d_gx_g,\]
 for some integer $d_g$.  We will prove that $d_g=1$,
which will of course prove Proposition 1.7. Set
$y_g:=f_g^*(x_g); f_g^*:H^2(M_g, \mathbb{Z}) \to H^2(J_g, \mathbb{Z}) $
 being the map in cohomology induced from $f_g$.

\begin{lemma}\label{L:nonvanish}
$y_g \neq 0$ and $r_g^*(y_{g+1})=y_g$ as elements of $H^2(J_g, \mathbb{Z}). $
\end{lemma}

\begin{proof}
There exists a unique universal line bundle $\mathcal{P}$, called the
{\it Poincar\'e bundle}
  on $\mathcal{C}_g \times \mathfrak{J}_{\mathcal{C}_g}$ such that, for each
$\mathfrak{L} \in \mathfrak{J}_{\mathcal{C}_g}$,  $\mathcal{P}$ restricts
 to the line bundle $\mathfrak{L}$  on $\mathcal{C}_g \times \mathfrak{L}$,
and   $\mathcal{P}$ restricted to $x_o\times \mathfrak{J}_{\mathcal{C}_g}$ is trivial for a fixed base point $x_o\in \mathcal{C}_g $ (cf. [ACGH, Chap. IV, $\S$2]).

Let $\mathcal{F}$ be the rank-$2$ vector bundle
$ \mathcal{P}\oplus \hat{\tau}^* (\mathcal{P})$ over the base space
$\mathcal{C}_g\times \mathfrak{J}_{\mathcal{C}_g},$ and think of
 $\mathcal{F}$ as a family of rank-$2$ bundles on $\mathcal{C}_g$ parametrized by $\mathfrak{J}_{\mathcal{C}_g},$ where $\hat{\tau}:\mathcal{C}_g\times \mathfrak{J}_{\mathcal{C}_g} \to \mathcal{C}_g\times \mathfrak{J}_{\mathcal{C}_g}$
is the involution $I\times \tau$.

By Drezet-Narasimhan [DN], we have
$x_g=c_1(\Theta_{V_2}(\mathcal{C}_g, SL(2)))$ for the standard representation $V_2$
of $SL(2)$.  Using the functoriality of Chern class,
\[ \xi^*_{\mathcal{C}_g}(x_g)=c_1(\mbox{Det}\,\mathcal{F}),\tag{1} \]
where $\mbox{Det}\,\mathcal{F}$ denotes the determinant line bundle over
$\mathfrak{J}_{\mathcal{C}_g}$ associated to the family  $\mathcal{F}$ (cf.
Definition 1.2). Recall that  the fiber of $\mbox{Det}\,\mathcal{F}$ at any
$\mathfrak{L} \in \mathfrak{J}_{\mathcal{C}_g}$ is
given by the expression
\[\begin{aligned}
{\mbox{Det}\,\mathcal{F}}_{|\mathfrak{L}}= &\wedge^{top}\bigl(H^0(\mathcal{C}_g, \mathfrak{L} \oplus \mathfrak{L}^{-1})^*\bigr) \otimes
       \wedge^{top} \bigl(H^1(\mathcal{C}_g, \mathfrak{L} \oplus
\mathfrak{L}^{-1})\bigr)\\
     =&\wedge^{top}\bigl(H^0(\mathcal{C}_g, \mathfrak{L})^* \oplus
      H^0(\mathcal{C}_g, \mathfrak{L}^{-1})^*\bigr) \otimes
      \wedge^{top}\bigl(H^1(\mathcal{C}_g, \mathfrak{L}) \oplus H^1(\mathcal{C}_g, \mathfrak{L}^{-1})\bigr) \\
     =&\wedge^{top}\bigl(H^0(\mathcal{C}_g, \mathfrak{L})^*\bigr) \otimes \wedge^{top}
\bigl(H^0(\mathcal{C}_g,
     \mathfrak{L}^{-1})^*\bigr) \otimes \wedge^{top}\bigl(H^1(\mathcal{C}_g,
\mathfrak{L})\bigr) \otimes
     \wedge^{top}(H^1(\mathcal{C}_g, \mathfrak{L}^{-1})) \\
    =&\bigl(\mbox{Det}\,\mathcal{P}\bigr)_{|\mathfrak{L}} \otimes
\bigl(\tau^*(\mbox{Det}\,\mathcal{P})\bigr)_{|\mathfrak{L}}.
\end{aligned}\tag{2} \]

Applying the  Grothendieck-Riemann-Roch theorem (cf. [F, Example 15.2.8]) for the projection
$\mathcal{C}_g \times  \mathfrak{J}_{\mathcal{C}_g}\leto{\pi} \mathfrak{J}_{\mathcal{C}_g}$ gives
\[ \ch(R\pi_* \mathcal{P})=\pi_*(\ch \mathcal{P}\cdot \Td T_{\pi}), \tag{3} \]
where $\ch$ is the Chern character and $\Td T_{\pi}$ denotes the Todd genus
of the relative tangent bundle of  $\mathcal{C}_g \times  \mathfrak{J}_{\mathcal{C}_g}$ along the fibers of $\pi$. By the definition of $\mbox{Det}\,\mathcal{P}$ and
$R\pi_* \mathcal{P}$,
\[
c_1(\mbox{Det}\,\mathcal{P})
                          =-\ch(R \pi_* \mathcal{P})_{[2]},\tag{4}
 \]
where, for a cohomology class $y$, $y_{[n]}$ denotes the component of $y$ in
$H^n$. Since  $\mathcal{P}$ restricted to $x_o\times \mathfrak{J}_{\mathcal{C}_g}$ is trivial and  for any
$\mathfrak{L} \in \mathfrak{J}_{\mathcal{C}_g}$,  $\mathcal{P}$ restricts
 to the line bundle $\mathfrak{L}$  on $\mathcal{C}_g \times \mathfrak{L}$
(with the trivial Chern class), we get
\[c_1(\mathcal{P})\in H^1(\mathcal{C}_g) \otimes
H^1(\mathfrak{J}_{\mathcal{C}_g}).\tag{5}\]
Thus, using (3)-(4),
\[\begin{aligned}
-c_1(\mbox{Det}\,\mathcal{P})&=\pi_*\left((\ch\,\mathcal{P}\cdot \Td T_{\pi})_{[4]} \right) \\
                         &=\pi_*\left( \frac{c_1(\mathcal{P})^2}{2}
                           + \frac{c_1(\mathcal{P}) \cdot c_1(T_{\pi})}{2}\right) \\
                         &=\pi_*\left(c_1(\mathcal{P})^2 \right)/2.
\end{aligned} \tag{6}\]
The last equality follows from (5), since  the cup product $c_1(\mathcal{P}) \cdot
c_1(T_{\pi})$ vanishes, $c_1(T_{\pi})$ being in $H^2(\mathcal{C}_g)\otimes
H^0(\mathfrak{J}_{\mathcal{C}_g})$.

Recall the presentation of $\pi_1(\mathcal{C}_g)$ given just above Theorem 1.5.
Then, $H_1(\mathcal{C}_g,  \mathbb{Z})=\oplus_{i=1}^{g} \mathbb{Z}a_i \hspace{0.03cm} \oplus \hspace{0.03cm}
\oplus_{i=1}^{g} \mathbb{Z}b_i$. Moreover,  the  $\Bbb Z$-module dual basis $\{a_i^*, b_i^* \}_{i=1}^{g}$ of $H^1(\mathcal{C}_g,  \mathbb{Z})=\Hom_{\Bbb Z}
(H_1(\mathcal{C}_g,  \mathbb{Z}), \Bbb Z)$
satisfies $a_i^* \cdot a_j^*=0=b_i^*\cdot b_j^*$, $a_i^*\cdot b_j^*=\delta_{ij}[\mathcal{C}_g]$, where $[\mathcal{C}_g]$ denotes the positive generator of
$H^2(\mathcal{C}_g, \mathbb{Z}).$

Having fixed a base point $x_o$ in $\mathcal{C}_g$, define the algebraic map
\[ \psi \; : \; \mathcal{C}_g \to  \mathfrak{J}_{\mathcal{C}_g}, \,\,
x \mapsto \mathcal{O}(x-x_o).\]

Of course, $ \mathfrak{J}_{\mathcal{C}_g}$ is canonically identified as
$H^1(\mathcal{C}_g, \mathcal{O}_{\mathcal{C}_g})/ H^1(\mathcal{C}_g,\Bbb Z)$.
Thus, as a real analytic space, we can identify
\[ \mathfrak{J}_{\mathcal{C}_g}\simeq  H^1(\mathcal{C}_g,\Bbb R)/
H^1(\mathcal{C}_g,\Bbb Z)\simeq  H^1(\mathcal{C}_g,\Bbb Z)\otimes_{\Bbb Z}
(\Bbb R/\Bbb Z)\simeq
\Hom_{\Bbb Z} \bigl(H_1(\mathcal{C}_g,\Bbb Z),  \Bbb R/\Bbb Z\bigr)=J_g \tag{7}\]
obtained from the $\Bbb R$-vector space isomorphism
\[ H^1(\mathcal{C}_g,\Bbb R)\simeq H^1(\mathcal{C}_g,
\mathcal{O}_{\mathcal{C}_g}),\]
induced from the inclusion $\Bbb R \subset \mathcal{O}_{\mathcal{C}_g}$, where
the last equality in (7) follows by using the basis $\{a_1,b_1,\ldots, a_g,b_g\}$ of $H_1(\mathcal{C}_g,\Bbb Z)$. The induced map, under the identification (7),
\[\psi_*:H_1(\mathcal{C}_g,\Bbb Z)\to H_1(\mathfrak{J}_{\mathcal{C}_g}, \Bbb Z)\simeq
H^1(\mathcal{C}_g,\Bbb Z)\]
is the Poincar\'e duality isomorphism. To see this, identify
\[\Hom_{\Bbb Z} \bigl(H_1(\mathcal{C}_g,\Bbb Z),  \Bbb R/\Bbb Z\bigr)\simeq
\Hom_{\Bbb Z} \bigl(H^1(\mathcal{C}_g,\Bbb Z),  \Bbb R/\Bbb Z\bigr)\tag{8}\]
using the Poincar\'e duality isomorphim: $ H_1(\mathcal{C}_g,\Bbb Z)\simeq
H^1(\mathcal{C}_g,\Bbb Z)$. Then, under the identifications (7)-(8), the map
 \[\psi :  \mathcal{C}_g \to \Hom_{\Bbb Z} \bigl(H^1(\mathcal{C}_g,\Bbb Z),  \Bbb R/\Bbb Z\bigr)\] can be described as
\[\psi(x)([\omega])=e^{2\pi i\int_{x_o}^x\omega},\]
for any closed $1$-form $\omega$ on $\mathcal{C}_g$ representing  the cohomology class $[\omega]\in H^1(\mathcal{C}_g,\Bbb Z)$ (cf. [M, Theorem 2.5]), where
$\int_{x_o}^x\omega$ denotes the integral of $\omega$ along any path in
$\mathcal{C}_g$ from $x_o$ to $x$.

Since
\[\psi_*: H_1(\mathcal{C}_g,\Bbb Z)\to  H_1(\mathfrak{J}_{\mathcal{C}_g}, \Bbb Z)\simeq H^1(\mathcal{C}_g,\Bbb Z)\]
is the Poincar\'e duality isomorphism, it is easy to see that the induced cohomology map
\[\psi^*:  H^1(\mathfrak{J}_{\mathcal{C}_g}, \Bbb Z)\simeq H_1(\mathcal{C}_g,\Bbb Z)
\to  H^1(\mathcal{C}_g, \Bbb Z) \]
is given by
\[\psi^*(a_i)=-b_i^*, \,\,\psi^*(b_i)=a_i^*\,\,\,\,\text{for all}\,\, 1\leq i\leq g.
\tag{9}\]
In particular, $\psi^*$ is  an isomorphism.
  Moreover, the isomorphism
does not depend on the choice of $x_o$.

Consider the map
\[ \mathcal{C}_g \times \mathcal{C}_g \leto{I \times \psi}
\mathcal{C}_g \times \mathfrak{J}_{\mathcal{C}_g}.\]
Let $\mathcal{P'}:=(I \times \psi)^*(\mathcal P)$.
Then, $\mathcal{P'}$ is the unique line bundle over $\mathcal{C}_g \times \mathcal{C}_g$
satisfying the following properties:

\[\mathcal{P'}|_{\mathcal{C}_g \times x}=\mathcal{O}(x-x_o) \;\; \mbox{and} \;\;
 \mathcal{P'}|_{x_o \times \mathcal{C}_g}\,\,\text{is trivial}.\]

Consider the following line bundle over $\mathcal{C}_g \times \mathcal{C}_g$:
\[ \mathcal{O}_{\mathcal{C}_g \times \mathcal{C}_g}(\triangle) \otimes
  (\mathcal{O}(-x_o) \boxtimes 1) \otimes (1\boxtimes \mathcal{O}(-x_o)), \]
where $\triangle$ denotes the diagonal in $\mathcal{C}_g \times \mathcal{C}_g$.
One sees that this bundle also satisfies the restriction properties mentioned above and hence
it must be isomorphic with $\mathcal{P'}$.  Consequently,
\[c_1(\mathcal{P'})=c_1(\mathcal{O}_{\mathcal{C}_g \times \mathcal{C}_g}(\triangle))+
  c_1(\mathcal{O}(-x_o) \boxtimes 1)+c_1(1\boxtimes \mathcal{O}(-x_o)) .\]

Using the definition of $\mathcal{P'}$ and the functoriality of the Chern
classes,
\[ c_1(\mathcal{P'})=c_1((I \times \psi)^*(\mathcal P))=(I \times \psi)^*c_1(\mathcal P). \tag{10}\]
By (5),  $c_1(\mathcal{P})\in H^1(\mathcal{C}_g)\otimes H^1( \mathfrak{J}_{\mathcal{C}_g})$, and hence
$c_1(\mathcal{P'})\in H^1(\mathcal{C}_g)\otimes H^1(\mathcal{C}_g)$.
Moreover,
\[ c_1(\mathcal{O}(-x_o) \boxtimes 1)+c_1(1\boxtimes \mathcal{O}(-x_o))\in
 H^2(\mathcal{C}_g)\otimes H^0(\mathcal{C}_g)\oplus  H^0(\mathcal{C}_g)\otimes H^2(\mathcal{C}_g).\]
Thus, $c_1(\mathcal{P'})$ is the component of
$c_1(\mathcal{O}_{\mathcal{C}_g \times \mathcal{C}_g}(\triangle))$ in $H^1(\mathcal{C}_g)\otimes H^1(\mathcal{C}_g)$.
Hence,  by Milnor-Stasheff [MS, Theorem 11.11],
\[c_1(\mathcal{P'})=- \ds \sum_{i=1}^{g} a_i^* \otimes b_i^* + \ds \sum_{i=1}^{g} b_i^*
                      \otimes a_i^*.\]
Therefore, by (10),
\[c_1(\mathcal{P})=- \ds \sum_{i=1}^{g} a_i^* \otimes {\psi^*}^{-1}(b_i^*) +
                     \ds \sum_{i=1}^{g} b_i^* \otimes {\psi^*}^{-1}(a_i^*),\]
and thus, by (6),
\[\begin{aligned}
c_1(\mbox{Det}\,\mathcal{P})&=-\frac{1}{2}\pi_*(c_1(\mathcal{P})^2) \\
                     &=-\frac{1}{2}\pi_* \Bigl(\bigl( -\ds \sum_{i=1}^{g} a_i^*
\otimes {\psi^*}^{-1}
                     (b_i^*)+ \ds \sum_{i=1}^{g} b_i^* \otimes {\psi^*}^{-1}
(a_i^*) \bigr)^2\Bigr) \\
                     &=-\frac{1}{2}\pi_* \left( \ds \sum_{i=1}^{g} a_i^* \cdot b_i^*  \otimes
                     {\psi^*}^{-1}(b_i^*)\cdot {\psi^*}^{-1}(a_i^*) +
                      \ds \sum_{i=1}^{g} b_i^* \cdot a_i^* \otimes {\psi^*}^{-1}(a_i^*)
                       \cdot {\psi^*}^{-1}(b_i^*) \right) \\
                    &=-\ds \sum_{i=1}^{g} {\psi^*}^{-1}(b_i^*) \cdot
{\psi^*}^{-1}(a_i^*)
                    \;\;\; \in H^2( \mathfrak{J}_{\mathcal{C}_g}, \mathbb{Z}).
\end{aligned} \]
Now, the involution $\tau$ of  $\mathfrak{J}_{\mathcal{C}_g}$ induces the map $-I$ on $H^1( \mathfrak{J}_{\mathcal{C}_g}, \mathbb{Z})$ (since, under the
identification $\bar{\theta}_{\mathcal{C}_g}:J_g \to  \mathfrak{J}_{\mathcal{C}_g}$, $\tau$ corresponds to the map $x\mapsto x^{-1}$ for $x\in J_g$). Therefore,
\[\tau^*(c_1(\mbox{Det}\,\mathcal{P}))=c_1(\mbox{Det}\,\mathcal{P}).\]
Hence, by the identities  (1)-(2),
\[ \begin{aligned}
 \xi^*_{\mathcal{C}_g}(x_g)=&c_1(\mbox{Det}\, \mathcal{F}) \\
     =&2c_1(\mbox{Det}\,\mathcal{P}) \\
     =&2 \ds \sum_{i=1}^{g} {\psi^*}^{-1}(a_i^*) \cdot {\psi^*}^{-1}(b_i^*),
\end{aligned} \tag{11}\]
which is clearly a nonvanishing class in $H^2(\mathfrak{J}_{\mathcal{C}_g}, \mathbb{Z})$. Moreover, for any $g\geq 2$, under the identification (7), the map
$r_{g-1}:J_{g-1} \to J_g$ corresponds to the map
$H_1(\mathcal{C}_{g}, \Bbb Z) \to H_1(\mathcal{C}_{g-1}, \Bbb Z), a_i\mapsto a_i,
b_i\mapsto b_i$ for $1\leq i\leq g-1,
a_g\mapsto 0, b_g \mapsto 0$. Thus, by (9) and (11), $ \xi^*_{\mathcal{C}_g}(x_g)$
 restricts, via $r_{g-1}^*$, to the class $
\xi^*_{\mathcal{C}_{g-1}}(x_{g-1})$ for any $g\geq 2$. But, by the commutative diagram (E), $\xi^*_{\mathcal{C}_g}(x_g)=y_g$.
 This proves Lemma~\ref{L:nonvanish}.
\end{proof}

\begin{proof}[Proof of Proposition 1.7]
By the above Lemma~\ref{L:nonvanish} and the commutative diagram (F), we see that
\[
f_g^*(d_gx_g)= f_g^*i_g^*(x_{g+1})=r_g^*(f_{g+1}^*(x_{g+1})),
\,\,\text{i.e.,}\,\,
d_gy_g= y_g. \]
Since the cohomology  of $J_g$ is torsion free and $y_g$ is a nonvanishing class,
we get $d_g=1$.  This concludes the proof of   Proposition~\ref{P:sl2case}.
\end{proof}


\newpage
\section*{References}
\[ \begin{array}{ll}
\mbox{[A]} & \mbox{Adams, F., {\it Lectures on Lie Groups}, W.A. Benjamin, Inc., 1969.}\\
\mbox{[ACGH]} & \mbox{Arbarello, E., Cornalba, M., Griffiths, P. and Harris, J., {\it Geometry}} \\
              & \mbox{{\it of Algebraic Curves}, Springer-Verlag, 1985.}\\
\mbox{[B]} & \mbox{Bourbaki, N., {\it Groupes et Alg\`ebres de Lie}, Chap. 4--6, Masson, Paris,}\\
            & \mbox{1981.}\\
\mbox{[BL]} & \mbox{Beauville, A. and Laszlo, Y., {\it Conformal blocks and generalized theta}}\\
            & \mbox{{\it functions}, Commun. Math. Phys. {\bf 164}, 385--419, (1994).}\\
\mbox{[BLS]} & \mbox{Beauville, A., Laszlo, Y. and Sorger, C., {\it The Picard group of the moduli}}\\
            & \mbox{{\it of G--bundles on a curve}, Compositio Math. {\bf 112}, 183--216, (1998).}\\
\mbox{[BR]} & \mbox {Beltrametti, M. and Robbiano, L., {\it Introduction to the theory of}}\\
            & \mbox {{\it weighted projective spaces}, Expo. Math. {\bf 4}, 111--162, (1986).}\\
\mbox{[BSh]} & \mbox{Bernshtein, I.N. and Shvartsman, O.V., {\it Chevalley's Theorem for complex}}\\
            & \mbox{{\it crystallographic Coxeter groups}, Funct. Anal. Appl. {\bf 12}, 308--310, (1978).}\\
\mbox{[DN]} & \mbox{Drezet, J.--M. and Narasimhan, M.S., {\it Groupe de Picard des vari\'et\'es de}}\\
            & \mbox{{\it modules de fibr\'es semi--stables sur les courbes alg\'ebriques}, Invent. Math.}\\
           & \mbox{{\bf 97}, 53--94, (1989).}\\
\mbox{[D]} &\mbox{Dynkin, E.B., {\it Semisimple subalgebras of semisimple Lie algebras}, Am. }\\
            & \mbox{Math. Soc. Transl. (Ser. II) {\bf 6}, 111--244, (1957).}\\
\mbox{[Fa]}& \mbox{Faltings, G., {\it A proof for the Verlinde formula}, J. Alg. Geom. {\bf 3}, 347--374,}\\
            & \mbox{(1994).}\\
\mbox{[F]} & \mbox{Fulton, W., {\it Intersection Theory}, Springer, 1998.}\\
\mbox{[FMW]} & \mbox{Friedman, R., Morgan, J.W. and Witten, E., {\it Principal $G$-bundles over }}\\
            & \mbox{{\it elliptic curves}, Mathematical Research Letters {\bf 5}, 97--118, (1998).}\\
\mbox{[KM]} & \mbox{Knudsen, F. and Mumford, D., {\it The projectivity of the moduli space of}} \\
& \mbox{{\it stable curves I: preliminaries on ``det'' and ``div''}, Math. Scand. {\bf 39}, 19--55,
(1976).}
\\
\mbox{[K$_1$]} & \mbox{Kumar, S., {\it Infinite Grassmannians and moduli spaces of $G$--bundles,} in:}\\
            & \mbox{{\it Vector Bundles on Curves - New Directions}, Lecture Notes in Mathematics }\\
            &\mbox{{\bf 1649}, Springer-Verlag, 1--49, (1997).}\\
\mbox{[K$_2$]} & \mbox{Kumar, S., {\it Kac-Moody Groups, their Flag Varieties and Representation}}\\
            & \mbox{{\it Theory}, Progress in Mathematics vol. {\bf 204}, Birkhauser, 2002.}\\
\mbox{[KN]} & \mbox{Kumar, S. and Narasimhan, M.S., {\it Picard group of the
moduli spaces of}}\\
            & \mbox{{\it G--bundles}, Math. Annalen {\bf 308}, 155--173, (1997).}
\\
\mbox{[KNR]} & \mbox{Kumar, S., Narasimhan, M.S., and Ramanathan, A., {\it Infinite Grassmannians}}\\
            &\mbox{{\it and moduli spaces of $G$--bundles}, Math. Annalen {\bf 300}, 41--75, (1994).}\\
\mbox{[L]} & \mbox{Lang, S., {\it Introduction to Arakelov Theory},
Springer, 1988.}\\
\mbox{[La]} & \mbox{Laszlo, Y., {\it About $G$-bundles over elliptic curves}, Ann. Inst. Fourier, Grenoble}\\
            & \mbox{ {\bf 48}, 413--424, (1998).}
\end{array}\]
\[\begin{array}{ll}

\mbox{[Lo]} & \mbox{Looijenga, E., \mbox{\it Root systems and elliptic curves,
Invent. Math. {\bf 38},}}\\
       &\mbox{17--32, (1976).}\\

\mbox{[LS]} & \mbox{Laszlo, Y. and Sorger, C., {\it The line bundles on the moduli of parabolic}}\\
        &\mbox{{\it G--bundles over curves and their sections}, Ann. Scient. \'Ec. Norm. Sup.}\\
        &\mbox{{\bf 30}, 499--525, (1997).}\\
\mbox{[M]} & \mbox{Milne, J.S., Chap. VII -- {\it Jacobian varieties,} in:
{\it Arithmetic Geometry}}\\
            & \mbox{ (Cornell, G., et al., eds.),
            Springer-Verlag, 167--212, (1986).}\\
\mbox{[MS]} & \mbox{Milnor, J.W. and Stasheff, J.D., {\it Characteristic Classes}, Annals of Mathematics}\\
        &\mbox{Studies vol. 76, Princeton University Press, 1974.}\\
\mbox{[N]} & \mbox{Narasimhan, R., {\it  Compact Riemann Surfaces},
Birkhauser Verlag, 1992.}\\
\mbox{[NS]} & \mbox{Narasimhan, M.S. and Seshadri, C.S., {\it Stable and unitary vector bundles on}}\\
        &\mbox{{\it a compact Riemann surface}, Ann. of Math. {\bf 82}, 540--567, (1965).}\\
\mbox{[RR]} & \mbox{Ramanan, S. and Ramanathan, A., {\it Some remarks on the instability flag},}\\
        &\mbox{T\^{o}hoku Math. J. {\bf 36}, 269--291, (1984).}\\
\mbox{[R$_1$]} & \mbox{Ramanathan, A., {\it Stable principal bundles on a compact Riemann surface-}}\\
            &\mbox{{\it Construction of moduli space}, Thesis, University of Bombay, 1976.}\\
\mbox{[R$_2$]} & \mbox{Ramanathan, A., {\it Stable principal bundles on a compact Riemann surface},}\\
            & \mbox{Math. Annalen {\bf 213}, 129--152, (1975).}\\
\mbox{[So]} & \mbox{Sorger, C., {\it On moduli of $G$-bundles on a curve for exceptional $G$,}}\\
        &\mbox{ Ann. Scient. \'Ec. Norm. Sup.}
        \mbox{\,{\bf 32}, 127--133, (1999).}\\
\mbox{[S]} & \mbox{Szenes, A., {\it The combinatorics of the Verlinde formula}, in: {\it Vector Bundles}}\\
            &\mbox{{\it in Algebraic Geometry} (Hitchin, N.J., et al., eds.), CUP {\bf 241}, (1995).}\\
\mbox{[T$_1$]} & \mbox{Teleman, C., {\it Lie algebra cohomology and the fusion rules}, Commun. Math.}\\
            &\mbox{Phys. {\bf 173}, 265--311, (1995).}\\
\mbox{[T$_2$]} & \mbox{Teleman, C., {\it Borel-Weil-Bott theory on the moduli stack of $G$-bundles over}}\\
            &\mbox{{\it a curve}, Invent. Math. {\bf 134}, 1--57, (1998).}\\
\mbox{[TUY]} & \mbox{Tsuchiya, A., Ueno, K. and Yamada, Y., {\it Conformal field theory on universal}}\\
            &\mbox{{\it family of stable curves with gauge symmetries}, Adv. Stud. Pure Math. {\bf 19},}\\
            &\mbox{ 459--565, (1989).}

\end{array} \]
\vskip5ex

\noindent
Address (of both the Authors):  Department of Mathematics,
University of North
Carolina,
Chapel Hill, NC 27599-3250, USA\\
shrawan$@$email.unc.edu\\
boysal$@$email.unc.edu
\end{document}